\documentclass[12pt]{article}
\usepackage[utf8]{inputenc}
\usepackage[english]{babel}
\usepackage[margin=1in]{geometry}
\usepackage{color}
\usepackage{enumerate}
\usepackage{amsmath, amsthm, amssymb, amsfonts}
\usepackage{bbm}
\usepackage{mathtools}
\usepackage[all,cmtip]{xy}
\usepackage{todonotes}
\usepackage{hyperref}
\allowdisplaybreaks
\newtheorem{theorem}{Theorem}[section]

\newtheorem{lemma}[theorem]{Lemma}

\newtheorem*{claim*}{Claim}

\newtheorem{proposition}[theorem]{Proposition}
\theoremstyle{definition}
\newtheorem{definition}[theorem]{Definition}

\newcommand{\Z}{\mathbb{Z}}
\newcommand{\F}{\mathbb{F}}

\newcommand{\Up}{\operatorname{Up}}

\newcommand{\syl}{\operatorname{Syl}}

\newcommand{\mbf}{\mathbf}
\newcommand{\sym}{\operatorname{Sym}}
\newcommand{\asym}{\operatorname{Antisym}}
\newcommand{\asymo}{\operatorname{Antisym0}}

\newcommand{\asymc}{\operatorname{antisym*}}

\newcommand{\divides}{\bigm|}

\newenvironment{theorem*}[2][Theorem]{\begin{trivlist}
\item[\hskip \labelsep {\bfseries #1}\hskip \labelsep {\bfseries #2}]}{\end{trivlist}}
\newenvironment{lemma*}[2][Lemma]{\begin{trivlist}
\item[\hskip \labelsep {\bfseries #1}\hskip \labelsep {\bfseries #2}]}{\end{trivlist}}
\newenvironment{corollary*}[2][Corollary]{\begin{trivlist}
\item[\hskip \labelsep {\bfseries #1}\hskip \labelsep {\bfseries #2}]}{\end{trivlist}}
\begin{document}
\title{\texorpdfstring{The Sylow $p$-subgroups of the finite simple groups}{The Sylow $p$-subgroups of the finite simple groups of classical Lie type}}
\author{Hannah Knight\thanks{This work was supported by NSF Award No. 2302822 as well as NSF Grant Nos. DMS-1811846 and DMS-1944862.}} 
\date{}
\maketitle

\maketitle

\vspace{-1cm}

\begin{abstract}

In this short article, we give a summary of the Sylow $p$-subgroups of the finite simple groups of classical Lie type. 
\end{abstract}



In \cite{Kni} and \cite{Kni2}, for $p \neq l$ primes, $q = p^r$, we have found the Sylow $p$-subgroups and the Sylow $l$-subgroups of $PSL_n(\F_q)$, $PSp(2n,q)$, and $P\Omega(2m,p^r)$, $\Omega(2m,2^r)$, and $\Omega(2m+1,p^r)$ as well as the Sylow $l$-subgroups of $PSU(n,p^r)$. In this summary, we also include the Sylow $p$-subgroups of $PSU(n,p^r)$.

\begin{lemma}[\cite{Kni2}, Lemma 2.18] Let $\sigma_{i}^{j}$ be the permutation which permutes the $i$th set of $l$ blocks of size $l^{j-1}$. Then 
$$\langle \{\sigma_i^j\}_{1 \leq j \leq \mu_l(n), 1 \leq i \leq \lfloor \frac{n}{l^{j}} \rfloor} \rangle \in \syl_l(S_n).$$ Let $P_l(S_n)$ denote this particular Sylow $l$-subgroup of $S_n$.
\end{lemma}

\begin{definition}
Let $v_l(n)$ denote the largest integer $i$ such that $l^i \divides n$. Let $\zeta_n$ denote a primitive $n$th root of unity. 
\end{definition}

\begin{definition} Define $\Up_{n}(\F_{p^r})$ to be the unitriangular $n \times n$ matrices over $\F_{p^r}$ under multiplication. (Unitriangular matrices are upper triangular matrices with $1$'s on the diagonal).\end{definition}

\begin{definition} 
The generalized quaternion groups are groups of order $4n$ with the following presentation:

$$Q_{4n} = \langle w, v : w^{n} = v^2, w^{2n} = 1, vwv^{-1} = w^{-1} \rangle.$$

\end{definition}

\section{\texorpdfstring{$PSL_n(\F_q)$}{PSLn(Fq)}}

\begin{proposition} Let $p$ be a primes, $n \geq 2$. Then for $P \in \syl_p(PSL_n(\F_q))$, 
$$P \cong \Up_n(\F_{p^r}).$$\end{proposition}

\begin{proposition} Let $p$ be a prime, $q = p^r$, $l$ a prime with $l \neq p$, and $s = v_l(q-1)$. Then for $P \in \syl_l(PSL_n(\F_q))$, 
$$P \cong \{(\mbf{b}, \tau) \in (\mu_{l^s})^n/ \{(x,\dots,x) : x^{n} = 1\} \rtimes P_l(S_n) : \prod_{i=1}^n b_i = \text{sgn}(\tau)\},$$
where the action of $P_l(S_n)$ on $(\mu_{l^s})^n$ is given by permuting the indices.
.\end{proposition}

\section{\texorpdfstring{$PSp(2n,p^r)$}{PSp(2n,pr)}}

\begin{definition} For any prime $p$, define $\sym(n,p^r)$ as the group of $n \times n$ symmetric matrices under addition (with entries from $\F_{p^r}$). \end{definition}

\begin{proposition} For any prime $p$, $P \in \syl_p(PSp(2n,p^r)$,
\begin{align*}
P &\cong \{\begin{pmatrix} A & 0_n\\
0_n & (A^{-1})^T\end{pmatrix}\begin{pmatrix} \text{Id}_n & B\\
0_n & \text{Id}_n\end{pmatrix} : A \in \Up_{n}(\F_{p^r}), B \in \sym(n,p^r)\}\\
&\cong \sym(n,p^r) \rtimes \Up_{n}(\F_{p^r}),\end{align*}
where the action of $A \in \Up_n(\F_{p^r})$ on $B \in \sym(n,p^r)$ is given by $A(B) = ABA^T$.
\end{proposition}

\begin{proposition}
Let $p$ be a prime and $l$ a prime with $l \neq 2,p$.  Let $d$ be the smallest positive integer such that $l \divides q^d - 1$. If $d$ is even then the Sylow $l$-subgroups of $PSp(2n,q)$ are isomorphic to Sylow $l$-subgroups of $GL_{2n}(\F_q)$. And if $d$ is odd, then the Sylow $l$-subroups of $PSp(2n,q)$ are isomorphic to Sylow $l$-subgroups of $GL_n(\F_q)$.
\end{proposition}

\begin{proposition} Let $p$ be a prime, $q = p^r$, and $l$ a prime with $l \neq p$. Let $d$ be the smallest positive integer such that $l \divides q^d-1$, $s = v_l(q^d-1)$, and $n_0 = \lfloor \frac{n}{d} \rfloor.$ Then for $P \in \syl_l(GL_n(\F_q))$, 
$$P \cong (\mu_{l^s})^{n_0} \rtimes P_l(S_{n_0}).$$
where the action of $P_l(S_{n_0})$ on $\mbf{b} \in (\mu_{l^s})^{n_0}$ is given by permuting the $b_i$.
\end{proposition}

\begin{proposition} Let $p \neq 2$ be prime. Then for $P \in \syl_2(PSp(2n,p^r)$
$$P \cong (Q_{2^s})^n/\langle (w^{2^{s-2}}, \dots w^{2^{s-2}}) \rangle \rtimes P_2(S_n).$$ 
where the action of $P_2(S_n)$ on $\mbf{a}$ is given by permuting the $a_i$.
\end{proposition}

\section{\texorpdfstring{$P\Omega(2m,p^r)$, $\Omega(2m,2^r)$, and $\Omega(2m+1,p^r)$}{POmega(2m,pr),Omega(2m,2r), and Omega(2m+1,pr)}}

\begin{definition} For any prime $p$, define $\asym(m,p^r)$ as the group of $m \times m$ anti-symmetric matrices under addition (with entries from $\F_{p^r}$).\end{definition}

\begin{definition}  For $p = 2$, define $\asymo(m,2^r) \subset \asym(m,2^r) = \sym(m,2^r)$ as the subgroup of symmetric/antisymmetric matrices with 0's on the diagonal. That is, $$\asymo(m,2^r) = \{B \in \sym(m,2^r) = \asym(m,2^r) : B_{i,i} = 0, \text{ } \forall i\}.$$\end{definition}

\begin{proposition} Let $p \neq 2$. Then for $P \in \syl_p(P\Omega^\epsilon(2m,p^r))$,
$$P \cong \asym(m,p^r) \rtimes \Up_{m}(\F_{p^r}),$$
where the action of $A \in \Up_n(\F_{p^r})$ on $B \in \asym(m,p^r)$ is given by $A(B) = ABA^T.$
\end{proposition}

\begin{proposition}
For $P \in \syl_2(\Omega^\epsilon(2m,2^r))$,
$$P \cong \asymo(m,2^r) \rtimes \Up_m(\F_{2^r}),$$ where the action of $A \in \Up_n(\F_{p^r})$ on $B \in \asymo(m,p^r)$ is given by
$A(B) = ABA^T.$ 
\end{proposition}

Note that $O^\epsilon(2m+1,2^r) \cong Sp(2m,2^r)$. So we only consider $p \neq 2$ for the odd orthogonal groups.
\begin{proposition}
 Let $p \neq 2$. Then for $P \in \syl_p(\Omega(2m+1,p^r))$,
$$P \cong \left((\F_{p^r}^+)^m \times \asym(m,p^r)\right) \rtimes \Up_{m}(\F_{p^r})),$$
where the action of $A \in \Up_m(\F_{p^r})$ on $B \in \asym(m,p^r)$ is given by $A(B) = ABA^T.$ and the action of $A \in \Up_m(\F_{p^r})$ on $\mbf{x} \in (\F_{p^r}^+)^m$ is given by $A(\mbf{x}) = \mbf{x}A^T.$
\end{proposition}

\begin{proposition}\label{On} Let $p$ be a prime, $q = p^r$, and $l$ a prime with $l \neq 2,p$.  Let $d$ be the smallest positive integer such that $l \divides q^d - 1$, and let $n_0 = \lfloor \frac{n}{d} \rfloor$. Then the Sylow $l$-subgroups of $P\Omega^\epsilon(n,q)$ are isomorphic to Sylow $l$-subgroups of
{\small \begin{align*}
\begin{cases} GL_m(\F_q), &n = 2m+1, \text{ d odd}\\
&\text{or } n = 2m, d \text{ odd}, \epsilon = + \\
GL_{m-1}(\F_q), &n = 2m, d \text{ odd}, \epsilon = - \\
GL_{2m}(\F_q),l, &n = 2m+1, \text{ d even}\\
&\text{or } n = 2m, \text{ d even}, n_0 \text{ even}, \epsilon = +\\
&\text{or } n = 2m, d \text{ even}, n_0 \text{ odd}, \epsilon = -\\
GL_{2m-2}(\F_q),l, &n = 2m, d \text{ even}, n_0 \text{ odd}, \epsilon = +\\
&\text{or } n = 2m, d \text{ even}, n_0 \text{ even}, \epsilon = -\\
\end{cases}
\end{align*}}
\end{proposition}

\begin{definition} 
The dihedral groups are groups of order $2n$ with the following presentation:
$$D_{2n} = \langle x,y : x^{n} = 1 = y^2, yxy = x^{-1} \rangle.$$
\end{definition}

\begin{proposition}\label{evenOsyl}  For $P \in \syl_2(O^+(2m,q))$,
\begin{align*}
P &\cong \begin{cases} (D_{2^{s+1}})^m \rtimes P_2(S_m), &q \equiv 1 \mod 4\\
&q \equiv 3 \mod 4, \text{ } m \text{ even}\\
\Z/2\Z \times \Z/2\Z \times ((D_{2^{s+1}})^{m-1}\rtimes P_2(S_{m-1})), &q \equiv 3 \mod 4, \text{ } m \text{ odd}\\
\end{cases}
\end{align*}
And for $P \in \syl_2(O^-(2m,q))$,
\begin{align*}
P &\cong \begin{cases} \Z/2\Z \times \Z/2\Z \times ((D_{2^{s+1}})^{m-1} \rtimes P_2(S_{m-1})), &q \equiv 1 \mod 4\\
&q \equiv 3 \mod 4, \text{ } m \text{ even}\\
(D_{2^{s+1}})^m \rtimes P_2(S_m), &q \equiv 3 \mod 4, \text{ } m \text{ odd}
\end{cases}
\end{align*}
\end{proposition}

\begin{proposition} For $P \in \syl_2(O(2m+1),2)$, $P \cong \Z/2\Z \times P'$ for 
$$P' \in \begin{cases} \syl_2(O^+(2m,q)), &q \equiv 1 \mod 4\\
&q \equiv 3 \mod 4, \text{ } m \text{ even}\\
\syl_2(O^-(2m,q)), &q \equiv 3 \mod 4, \text{ } m \text{ odd}
\end{cases}.$$
\end{proposition}

\section{\texorpdfstring{$PSU(n,p^r)$}{PSU(2n,pr)} }

 \begin{definition} Let $\asymc(m,p^{2r}) = \{B \in M_{n\times n}(\F_{p^{2r}}) : B^T = -\overline{B}\}$ under addition. \end{definition} 
 
The kernel of the natural homomorphism  $U(n,p^{2r}) \to PSU(n,p^{2r})$ has order prime to $p$, so it maps the Sylow $p$-subgroups of $U(n,p^{2r})$ isomorphically onto Sylow $p$-subgroups of $PSU(n,p^{2r})$, so it suffices to consider the Sylow $p$-subgroups of $U(n,p^{2r})$. It is a straightforward calculation to check the following propositions.
 \begin{proposition} For $n = 2m$, any $p$, $P \in \syl_p(U(2m,p^{2r}))$,
 \begin{align*}
 P &\cong \{\begin{pmatrix} A & 0_m\\
0_m & (\overline{A^{-1}})^T\end{pmatrix}\begin{pmatrix} \text{Id}_m & B\\
0_m & \text{Id}_m\end{pmatrix} : A \in \Up_{m}(\F_{p^{2r}}), B \in \asymc(m,p^{2r})\}\\
&\cong \asymc(m,p^{2r}) \rtimes \Up_{m}(\F_{p^{2r}}),
\end{align*}
where the action of $A \in \Up_m(\F_{p^{2r}})$ on $B \in \asymc(m,p^{2r})$ is given by $A(B) = AB\overline{A}^T.$
\end{proposition}

 \begin{proposition} Let $$
S = \{(\mbf{y},B) : \mbf{y} \in (\F_{p^{2r}})^m, B \in M_{m \times m}(\F_{p^{2r}}), B^T + \overline{B} = -\mbf{y}^T\overline{\mbf{y}}\},$$ where the multiplication is given by $(\mbf{y},B)(\mbf{y}',B') = (\mbf{y}+\mbf{y'}, B+B' - \overline{\mbf{y}}^T\mbf{y}')$. 

Then for $P \in \syl_p(U(2m+1,p^{2r}))$, 
$$P \cong S \rtimes \Up_m(\F_{p^{2r}}),$$
where the action of $\mbf{A} \in \Up_m(\F_{p^{2r}})$ on $(\mbf{y},B) \in S$ is given by 
$$A(\mbf{y},B) = (\mbf{y}\overline{A}^T, AB\overline{A}^T).$$
 \end{proposition}

 \begin{proposition}[\cite{Kni2}, Sections 9.4 and 9.3] For $l \neq p$, $q = p^2$, the Sylow $l$-subgroups of $PSU(n,q^2)$ are isomorphic to Sylow $l$-subgroups of
$$\begin{cases}
PSL_n(\F_{q^2}), &l \divides n, l \divides q+1\\
SL_n(\F_{q^2}), &l \nmid n, l \divides q+1\\
U(n,q^2), &l \nmid n, l \nmid q+1
\end{cases}.$$
\end{proposition}

\begin{proposition}[\cite{Kni2}, Theorem 9.7] Let $l \neq p$ be primes and $q = p^r$. Let $d$ be the smallest positive integer such that $l \divides q^d-1$. Then Sylow $l$-subgroups of $U(n,q^2)$ are isomorphic to Sylow $l$-subgroups of 
$$\begin{cases} GL_n(\F_{q^2}), &d = 2 \mod 4\\
GL_{\lfloor \frac{n}{2}\rfloor}(\F_{q^2}), &d \neq 2 \mod 4 \end{cases}.$$
\end{proposition}

\begin{proposition}\label{GLsyl} Let $l \neq p$ be primes and $q = p^r$. Let $d$ be the smallest positive integer such that $l \divides q^d - 1$, $s = v_l(q^d-1)$, and $n_0 = \lfloor \frac{n}{d} \rfloor$.  For $P \in \syl_l(GL_n(\F_q))$,  
$$P \cong (\mu_{l^s})^{n_0} \rtimes P_l(S_{n_0}).$$
\end{proposition}

\begin{lemma} Let $l \neq p$ be primes and $q = p^r$. For $P \in \syl_l(SL_n(\F_q))$, 
$$P \cong \{(\mbf{b}, \tau) \in (\mu_{l^s})^n \rtimes P_l(S_n) : \prod_{i=1}^n b_i = \text{sgn}(\tau)\},$$
where the action of $P_l(S_n)$ on $\mbf{b} \in (\mu_{l^s})^n$ is given by permuting the $b_i$.
\end{lemma}

\bibliographystyle{plain}
\bibliography{references}

\end{document}